\documentclass{article}
\usepackage{graphicx} % Required for inserting images
\usepackage[utf8]{inputenc} % Required for inputting international characters
\usepackage{comment}
\usepackage{amsmath}
\usepackage{amssymb}
\usepackage{mathtools}
\usepackage{amsthm}
\usepackage{thmtools}
\usepackage{thm-restate}
\usepackage{newtxmath}
\usepackage{stmaryrd} 
\usepackage[mathscr]{euscript}
 \let\mathscr\relax
\usepackage[scr]{rsfso}
\usepackage{comment}
\usepackage[colorlinks=true, allcolors=blue]{hyperref}
\usepackage{url}
\usepackage{cleveref}
\usepackage[notintoc]{nomencl}
\usepackage{enumitem}
\counterwithout{footnote}{section}
\newcommand{\olsi}[1]{\,\overline{\!{#1}}} % overline short italic
\newcommand{\closedb}{\olsi{B}}

\newcommand{\diam}{\normalfont{\text{diam}}}

\newcommand{\image}{\normalfont{\text{im}}}

\newcommand{\len}{\normalfont{\text{len}}}

\newcommand{\Dya}{{\{0,1\}}^{\Naturals}}

\newcommand{\Naturals}{\mathbb{N}}

\numberwithin{equation}{section}
\theoremstyle{plain}
\newtheorem{proposition}{Proposition}[section]
\newtheorem{lemma}[proposition]{Lemma}
\newtheorem{theorem}[proposition]{Theorem}
\newtheorem{corollary}[proposition]{Corollary}

\theoremstyle{definition}
\newtheorem{definition}[proposition]{Definition}
\newtheorem{notation}[proposition]{Notation}

\theoremstyle{remark}
\newtheorem{remark}[proposition]{Remark}

\newtheorem{example}[proposition]{Example}

\title{Ultrametrizable spaces are homeomorphic to clade spaces of pruned trees}
\author{Itamar Bella\"iche}
\date{July 2024}

\begin{document}

\maketitle
\begin{abstract}
    This paper demonstrates that every ultrametric space is homeomorphic to a clade space of a pruned tree, \emph{i.e.}, a subspace of a tree’s canopy. Furthermore, it characterizes several topological properties of ultrametrizable spaces through the features of their representing trees. This approach suggests that topological properties of ultrametrizable spaces should be studies via the study of naturally ordered pruned trees.
\end{abstract}

\section{Introduction}
This paper expands upon the third chapter of my master thesis \cite{bellaiche2023}, exploring the topological relationships between ultrametric spaces and pruned trees of depth $\omega$. While these connections have been previously examined from a categorical perspective for completely ultrametrizable spaces \cite{brian2015completely}, this work strengthens this analysis and extends it to a broader class of spaces.

The central result of this paper, presented in Corollary \ref{HomeomorphismCorollary}, establishes a topological equivalence between ultrametrizable spaces and \emph{clade} spaces of pruned trees (see Definition \ref{CladeCanopyDefinition}). This equivalence provides a powerful framework for understanding the topological properties of ultrametric spaces through the lens of pruned trees. Consequently, when investigating the topology of ultrametrizable spaces, it suffices to consider the special case of pruned trees and their associated clade spaces.

The paper is structured as follows: Section \ref{Preliminaries} presents some notations, definitions, and basic results about ultrametric spaces. Section \ref{Sequentially} present a special case of ultrametric spaces, called \emph{sequentially-descending-to-zero} ultrametric spaces, which will be useful for the main proofs. Finally, Section \ref{Trees} presents the realm of pruned trees of depth $\omega$, states the definition of clade spaces, proves the main result of the paper, and finishes by the study of the manifestation of some topological properties of ultrametrizable spaces in the structure of their representing trees. 

\section{Preliminaries}
\label{Preliminaries}
Given a topological space $(X,\mathcal{O})$, denote by $\mathcal{B}_{(X,\mathcal{O})}$ the Borel tribe\footnote{This paper uses the term ``tribe'' instead of ``Boolean $\sigma$-algebra'', as is done in the majority of francophone mathematical texts.} generated by the topology $\mathcal{O}$.
Given a metric space $(X,d)$, we define and denote the open ball of radius $r\in(0,\infty)$ centered at $x\in X$ by
\begin{equation}
	B(x,r)=\{y\in X\:|\:d(x,y)<r\}
\end{equation}
and the closed ball of radius $r\in(0,\infty)$ centered at $x\in X$ by
\begin{equation}
	\closedb(x,r)=\{y\in X\:|\:d(x,y)\leq r\}.
\end{equation}
The topology on $X$ induced by the metric $d$ is the topology generated by the basis consisting of all open balls in the space, and is denoted by $\mathcal{O}_{(X,d)}$. A topological space is said to be ``metrizable'' if there exists a metric that induces its topology.
Denote also $\mathcal{B}_{(X,d)}=\mathcal{B}_{\mathcal{O}_{(X,d)}}$.
The distance of an element $x\in X$ from a non-empty set $A\subseteq X$ is defined and denoted by $d(x,A)=\inf\{d(x,y)\:|\:y\in A\}$, and the distance between two non-empty sets $A,B\subseteq X$ is defined and denoted by $d(A,B)=\inf\{d(x,B)\:|\:x\in A\}$.
\begin{definition}
	\label{SeparatedSetDefinition}
	Let $(X,d)$ be a metric set, $\emptyset\neq A\subseteq X$, and $r\in[0,\infty)$. Then $A$ is said to be \emph{$r$-separated} if for every $x,y\in A$ we have $d(x,y)>r$.
\end{definition}
\begin{comment}
	\begin{proposition}[open set is $F_\sigma$]
		\label{FsigmaProposition}
		Let $(X,d)$ be a metric space. Then every open set (according to the induced topology) is $F_\sigma$, \emph{i.e.}, a union of countably many closed sets.
	\end{proposition}
	\begin{proof}
		Let $U\in\mathcal{O}_{(X,d)}$ be an open set according to the topology induced by the metric $d$. Then $U^c=X\setminus U$ is closed.\\
		For every $k\in\Naturals\setminus\{0\}$ denote the open set
		\begin{equation*}
			G_k=\left\{x\in X\:\middle|\:d(x,U^c)<\tfrac{1}{k}\right\}=\bigcup_{x\in U^c}B\left(x,\tfrac{1}{k}\right)
		\end{equation*}
		Clearly we have
		\begin{equation*}
			U^c=\bigcap_{k=1}^\infty G_k
		\end{equation*}
		and thus $U^c$ is $G_\delta$, \emph{i.e.}, an intersection of countably many open sets. As a consequence, by De Morgan's law, the open set $U$ is $F_\sigma$:
		\begin{equation*}
			U=X\setminus U^c=X\mathbin{\big\backslash}\bigcap_{k=1}^\infty G_k=\bigcup_{k=1}^\infty X\setminus G_k
		\end{equation*}
		and $X\setminus G_k$ is closed for every $k\in\Naturals\setminus\{0\}$.
	\end{proof}
\end{comment}
\subsection{Doubling metric space}
We finish this section by presenting the definition of a \emph{doubling} metric space. 
\begin{definition}
	\label{DoublingDefinition}
	Let $(X,d)$ be a metric space. It is a \emph{doubling} metric space if there exists a \emph{doubling constant} $D\in\Naturals\setminus\{0\}$ such that for every radius $r\in(0,\infty)$ and every point $x\in X$, the closed ball $\closedb(x,r)$ can be covered by at most $D$ closed balls of radius $\tfrac{r}{2}$.
\end{definition}
\begin{remark}
	Every doubling metric space is separable.
\end{remark}

\subsection{Ultrametric spaces}
\label{UltrametricSubsection}
This subsection presents some useful results that are true for general ultrametric spaces. An ultrametric is a metric that satisfies a stronger condition than the triangle inequality, which is known as the ``ultrametric inequality'' or the ``strong triangle inequality''. Here is its formal definition:
\begin{definition}[ultrametric]
	An ultrametric over $X$ is a non-negative function $d\colon X\times X\longrightarrow[0,\infty)$ satisfying:
	\begin{enumerate}
		\item (symmetry) $d(x,y)=d(y,x)$ for every $x,y\in X$
		\item (identity of indiscernibles) $d(x,y)=0$ if and only if $x=y\in X$
		\item (strong triangle inequality)  $d(x,y)\leq\max\{d(x,z),d(z,y)\}$ for every $x,y,z\in X$
	\end{enumerate}
\end{definition}
A pair $(X,d)$, where $d$ is an ultrametric over $X$, constitutes an ultrametric space. Since the ultrametric inequality is stronger than the triangle inequality, every ultrametric is a metric, and the definitions of open and closed balls, the induce topology, the Borel tribe generated by that topology and other metric notions remain the same. A topological space is said to be ``ultrametrizable'' if there exists an ultrametric that induces its topology.
\begin{proposition}[isosceles triangle with a short base]
	\label{isoscelesProposition}
	Let $(X,d)$ be a metric space. Then $d$ is an ultrametric on $X$ if and only if every three points in $(X,d)$ constitute an isosceles triangle with a short base, \emph{i.e.}, for every $x,y,z\in X$, if $d(x,y)<d(x,z)$ then $d(y,z)=d(x,z)$.
\end{proposition}
\begin{proof}
    Let $(X,d)$ be a metric space.\\
    Suppose that $d$ is an ultrametric on $X$.\\
    Let $x,y,z\in X$ such that $d(x,y)<d(x,z)$.\\
    Suppose by contradiction that $d(y,z)\neq d(x,z)$. Then by the strong triangle inequality,
    \begin{equation*}
        d(y,z)\leq\max\{d(x,y),d(x,z)\}=d(x,z).
    \end{equation*}
    Therefore, $d(y,z)<d(x,z)$. But then, invoking again the strong triangle inequality we get that
    \begin{equation*}
        d(x,z)\leq\max\{d(x,y),d(y,z)\}<d(x,z),
    \end{equation*}
    a contradiction. Thus $d(y,z)=d(x,z)$.\\
    Suppose now that for every $x,y,z\in X$, if
    \begin{equation*}
        d(x,y)<d(x,z),
    \end{equation*}
    then
    \begin{equation*}
        d(y,z)=d(x,z).
    \end{equation*}
    Let $x,y,z\in X$.\\
    Suppose by contradiction and w.l.o.g. that
    \begin{equation*}
        d(x,z)>\max\{d(x,y),d(y,z)\}.
    \end{equation*}
    Then $d(x,z)>d(x,y)$, and by the assumption on $(X,d)$, we have $d(y,z)=d(x,z)$, in contradiction to the assumption that also $d(x,z)>d(y,z)$. Thus
    \begin{equation*}
        d(x,z)\leq\max\{d(x,y),d(y,z)\},
    \end{equation*}
    and $d$ is indeed an ultrametric on $X$.
\end{proof}

\subsection{Ultrametric balls}
\label{UltrametricBallsSubsection}
Open and closed balls in ultrametric spaces have some basic special properties:
\begin{proposition}[ultrametric balls]
    \label{UltrametricBallsBasics}
    Let $(X,d)$ be an ultrametric space, $x\in X$ an element in it and $r\in(0,\infty)$ a radius. Then
    \begin{enumerate}
		\item\label{diamItem} $\diam B(x,r)\leq\diam\closedb(x,r)\leq r$.
		\item\label{closedDiamItem} $\closedb(x,r)=\closedb\left(x,\diam\closedb(x,r)\right)$.
		\item\label{egocentricity} If $y\in B(x,r)$ then $B(y,r)=B(x,r)$. The same holds also for closed balls. This means that all points in a ball are its centers. This property is sometimes called \emph{egocentricity}.
		\item\label{inclusionBalls} If $q\leq r$ and $y\in X$ such that $B(x,r)\cap B(y,q)\neq\emptyset$,  then $B(y,q)\subseteq B(x,r)$ and $B(x,r)= B(y,r)$. The same holds also for closed balls.
		\item\label{clopenItem} The open ball $B(x,r)$ is a closed set and the closed ball $\closedb(x,r)$ is an open set.
		\item\label{CompactClopen} If $(X,d)$ is also compact, then the non-empty open ball $B(x,r)$ is also a closed ball, and the closed ball $\closedb(x,r)$ is also an open ball.
	\end{enumerate}
\end{proposition}
\begin{proof}
	We prove these properties by their order.
	\begin{enumerate} 
		\item By definition,
        \begin{equation*}
            \begin{split}
                \diam B(x,r)&=\sup\{d(y,z)\:|\:y,z\in B(x,r)\}\\
                &\leq\sup\{d(y,z)\:|\:y,z\in\closedb(x,r)\}=\diam\closedb(x,r).\\
            \end{split}
        \end{equation*}
        By the ultrametric inequality, for every $y,z\in\closedb(x,r)$ we have
        \begin{equation*}
            d(y,z)\leq\max\{d(y,x),d(x,z)\}\leq r
        \end{equation*}
        and thus $\diam\closedb(x,r)\leq r$.
		\item From point \ref{diamItem} we have that
        \begin{equation*}
            \closedb\left(x,\diam\closedb(x,r)\right)\subseteq\closedb(x,r).
        \end{equation*}
        Let $y\in\closedb(x,r)$. Then
        \begin{equation*}
            d(x,y)\leq\sup\{d(z,w)\:|\: z,w\in\closedb(x,r)\}=\diam\closedb(x,r)
        \end{equation*}
        and thus
        \begin{equation*}
            y\in\closedb\left(x,\diam\closedb(x,r)\right).
        \end{equation*}
        Therefore,
        \begin{equation*}
            \closedb\left(x,\diam\closedb(x,r)\right)=\closedb(x,r).
        \end{equation*}
		\item Suppose $y\in B(x,r)$ and let $z\in B(x,r)$. Then by the ultrametric inequality,
        \begin{equation*}
            d(y,z)\leq\max\{d(x,y),d(x,z)\}<r,
        \end{equation*}
        and therefore $z\in B(y,r)$. Thus
        \begin{equation*}
            B(x,r)\subseteq B(y,r).
        \end{equation*}
        By the same reasoning we get that
        \begin{equation*}
            B(x,r)\supseteq B(y,r),
        \end{equation*}
        and thus $B(x,r)=B(y,r)$. The proof for the closed ball case is analogous.
		\item Suppose that $B(x,r)\cap B(y,q)\neq\emptyset$, and let $z\in B(x,r)\cap B(y,q)$. Then by the ultrametric inequality,
        \begin{equation*}
            d(x,y)\leq\max\{d(x,z),d(y,z)\}<\max\{r,q\}.
        \end{equation*}
        If $q\leq r$, then $\max\{r,q\}=r$ and $d(x,y)<r$. Thus $y\in B(x,r)$ and so by point \ref{egocentricity} we get $B(y,q)\subseteq B(y,r)=B(x,r)$. The proof for the closed balls case in analogous.
		\item We begin by showing that the open ball $B(x,r)$ is a closed set.\\
		Denote $B^c=X\setminus B(x,r)$ and let $y\in B^c$. Thus $d(x,y)\geq r$.\\
		Let $z\in B(y,r)$. Suppose by contradiction that $d(x,z)<r$. Then by the ultrametric inequality
        \begin{equation*}
            d(x,y)\leq\max\{d(x,z),d(z,y)\}<r,
        \end{equation*}
        a contradiction. Thus $d(x,z)\geq r$ and $z\in B^c$. We get that $B(y,r)\subseteq B^c$. Since this holds for every $y\in B^c$, we get that
		\begin{equation}
			\bigcup_{y\in B^c}B(y,r)=B^c
		\end{equation}
		and thus $B^c$ is open as a union of open balls. Hence $B(x,r)$ is closed in the topology induced by the ultrametric and $\overline{B(x,r)}=B(x,r)$.\\
		We now show that $\closedb(x,r)$ is an open set.\\
		Let $y\in\closedb(x,r)$. By point \ref{egocentricity} of the current proposition
        \begin{equation*}
            B(y,r)=B(x,r)\subseteq\closedb(x,r).
        \end{equation*}
        Thus
		\begin{equation}
			\bigcup_{y\in\closedb(x,r)}B(y,r)=\closedb(x,r)
		\end{equation}
		and the closed ball $\closedb(x,r)$ is thus open as it is a union of open balls.
		\item Suppose $(X,d)$ is also compact.\\
        Then we must have $\diam B(x,r)<r$:\\
		Suppose by contradiction that $\diam B(x,r)=r$.\\
		By the ultrametric inequality we get that
		\begin{equation*}
            \begin{split}
                &\sup\{d(z,y)\:|\:y,z\in B(x,r)\}=\sup\{d(x,y)\:|\:y\in B(x,r)\}\\
                &=\diam B(x,r)=r.\\
            \end{split}
		\end{equation*}
		Therefore, for every $k\in\Naturals\setminus\{0\}$, there exists some $z_k\in B(x,r)$ such that $d(x,z_k)>r-\tfrac{1}{k}$. Since the space $(X,d)$ is compact, there must exist a subsequence $(z_{k_j})_{j=1}^\infty$ that admits a limit, denote it by $z$. Since $B(x,r)$ is closed (according to point \ref{clopenItem}), we must have
        \begin{equation*}
            z=\lim_{j\to\infty}z_{k_j}\in B(x,r),
        \end{equation*}
        and thus $d(x,z)<r$. However, by the continuity of the distance function, we must have
        \begin{equation*}
            d(x,z)=\lim_{j\to\infty}d(x,z_{k_j})=r,
        \end{equation*}
        a contradiction.\\
		Hence $\diam B(x,r)<r$ and thus $B(x,r)\supseteq\closedb\bigl(x,\diam B(x,r)\bigr)$.\\
		Additionally, by point \ref{closedDiamItem} we have
		\begin{equation*}
			B(x,r)\subseteq\closedb(x,r)=\closedb\left(x,\diam B(x,r)\right)
		\end{equation*}
		and we get that
		\begin{equation*}
			B(x,r)=\closedb\left(x,\diam B(x,r)\right)
		\end{equation*}
		and it is indeed a closed ball.\\
		We now show that $\closedb(x,r)$ is also an open ball.\\
		If $\diam\closedb(x,r)<r$ then by point \ref{closedDiamItem} we get that
		\begin{equation*}
			B(x,r)\subseteq\closedb(x,r)=\closedb(x,\diam\closedb(x,r))\subseteq B(x,r)
		\end{equation*}
		and thus $\closedb(x,r)=B(x,r)$ and it is clearly also an open ball.
		Else, $\diam\closedb(x,r)=r$ and there exists a point $y\in X$ such that $d(x,y)=r$. Then $y\in B(x,q)$ and $\diam B(x,q)\geq r$ for every $q\in(r,\infty)$.\\
		Suppose by contradiction that $\diam B(x,q)>r$ for every $q\in(r,\infty)$.\\
		As we saw, since $(X,d)$ is compact, we must have $\diam B(x,q)<q$ for every $q\in(r,\infty)$. By the sandwich rule, we get that
		\begin{equation*}
			\lim_{q\to r^-}\diam B(x,q)=r.
		\end{equation*}
		For every $k\in\Naturals\setminus\{0\}$ let $z_k\in B(x,r+\tfrac{1}{k})$ such that
		\begin{equation*}
			d(x,z_k)=\diam B(x,r+\tfrac{1}{k})\in(r,r+\tfrac{1}{k}).
		\end{equation*}
		By the compactness of $(X,d)$, there exists a subsequence $(z_{k_j})_{j=1}^\infty$ that admits a limit, denote it by $z$. On one hand, by the definition of a limit, we get that $\lim d(z_{k_j},z)=0$. On the other hand, by the continuity of the ultrametric, we get that
		\begin{equation*}
			d(x,z)=\lim_{j\to\infty}d(x,z_{k_j})=\lim_{j\to\infty}\diam B(x,r+\tfrac{1}{k_j})=r.
		\end{equation*}
		However, by Proposition \ref{isoscelesProposition}. we also get that
		\begin{equation*}
			\begin{split}
				d(z,z_{k_j})&=\max\{d(x,z),d(x,z_{k_j})\}=\max\{ r,\diam B(x,r+\tfrac{1}{k_j})\}=\\
				&=\diam B(x,r+\tfrac{1}{k_j})
			\end{split}
		\end{equation*}
		for every $k\in\Naturals\setminus\{0\}$. Invoking again the continuity of the ultrametric, we get that
		\begin{equation*}
			0=d(z,z)=d(z,\lim_{j\to\infty}z_{k_j})=\lim_{j\to\infty}d(z,z_{k_j})=\lim_{j\to\infty}\diam B(x,r+\tfrac{1}{k_j})=r
		\end{equation*}	
		a contradiction.\\
		Thus there exists a number $q^*\in(r,\infty)$ such that $\diam B(x,q)=r$ for every $q\in(r,q^*]$. Thus $\closedb(x,r)=B(x,q^*)$ and it is indeed an open ball.
	\end{enumerate}
\end{proof}
\begin{remark}
	As a consequence of point \ref{clopenItem} of Proposition \ref{UltrametricBallsBasics}, we get that the topology induced by an ultrametric has a clopen basis. Thus, an ultrametrizable space has topological dimension zero and is totally disconnected, see \cite{brian2015completely}.
\end{remark}

\section{Sequentially-descending-to-zero ultrametrics}
\label{Sequentially}
\begin{comment}
	The term ``sequentially-desending-to-zero ultrametric'' is quite long and cumbersome. A better nomenclature is thus clearly needed. Considering proposition \ref{universalProposition}, the name ``universal ultrametric'' might be in place. The next section also suggest using the term ``genealogic'' or ``dendronic''\footnote{$\delta\epsilon\nu\delta\rho o\nu$ (``Dendron'') means ``tree'' in ancient greek.}  ultrametric.
\end{comment}
\begin{definition}[sequentially-descending-to-zero set]
	Let $A\subset[0,\infty)$ be a set. The set $A$ is said to be \emph{sequentially-descending-to-zero} if there exists a sequence of non-negative real numbers $(a_k)_{k\in\Naturals}\in[0,\infty)^\Naturals$ that converges to zero and such that
	\begin{equation}
		A=\{0\}\cup\{a_k \:|\: k\in\Naturals\}.
	\end{equation}
\end{definition}
Note that if $A$ is sequentially-descending-to-zero, then $|A|\leq\aleph_0$.
\begin{proposition}
	\label{SDTZsetProperty}
	Let $A\subset[0,\infty)$. Then $A$ is sequentially-descending-to-zero if and only if $A$ is bounded, $\min A=0$, and for every $r\in(0,\infty)$ there is an $\varepsilon_r>0$ such that $(r-\varepsilon_r,r+\varepsilon_r)\setminus\{r\}\cap A=\emptyset$.
\end{proposition}
\begin{proof}
	Let $A\subset[0,\infty)$.\\
	Suppose $A$ is sequentially-descending-to-zero. Let $r\in(0,\infty)$.\\
	Let $(a_k)_{k\in\Naturals}$ be w.l.o.g a weakly monotonically decreasing\footnote{one can transform $(a_k)_{k\in\Naturals}$ into a weakly monotonically decreasing sequence $(b_k)_{k\in\Naturals}$ by defining $b_k=\min\{a_j \:|\:0\leq j\leq k\}$ for every $k\in\Naturals$.} sequence that converges to zero such that $A=\{0\}\cup\{a_k \:|\: k\in\Naturals\}$. Let $K_r\in\Naturals$ such that $a_k<\frac{r}{2}$ for every $k\geq K_r$ and denote
	\begin{equation*}
		\varepsilon_r = \min\left\{|a_j-r|\:\middle|\: j\in\{0,1,...,K_r\}\text{ and }a_j\neq r\right\}
	\end{equation*}
	then indeed $(r-\varepsilon_r,r+\varepsilon_r)\setminus\{r\}\cap A=\emptyset$. Since $(a_k)_{k\in\Naturals}$ converges to zero, it is bounded and thus also the set $A$.\\
	Suppose now that $A$ is bounded, that $0\in A$, and that for every $r\in(0,\infty)$ there is an $\varepsilon_r>0$ such that $(r-\varepsilon_r,r+\varepsilon_r)\setminus\{r\}\cap A=\emptyset$. We show that $A$ is sequentially-descending-to-zero:\\
	Denote $R=\sup A$. Since there exists an $\varepsilon_R>0$ such that
    \begin{equation*}
        (R-\varepsilon_R,R+\varepsilon_R)\setminus\{R\}\cap A=\emptyset,
    \end{equation*}
    the supremum $R$ must be a maximum of $A$, \emph{i.e.}, $R\in A$. Therefore define $a_0=R$. By the same reasoning, we can continue inductively: for every $k\in\Naturals\setminus\{0\}$, denote
	\begin{equation*}
		a_k=\max\left((A\setminus\{a_j \:|\: 0\leq j< k\})\cup\{0\}\right).
	\end{equation*}
	The sequence $(a_k)_{k\in\Naturals}$ is clearly weakly monotonically decreasing and consists of elements from $A$. Thus it must converge to a limit, denote it by $a\geq0$.  Suppose by contradiction that $a>0$. Then, if there is $k_a\in\Naturals$ such that $a_{k_a}=a$, by the inductive construction of $(a_k)_{k\in\Naturals}$ we must have $a_{k_a+1}=0$, in contradiction to the assumption that $\lim_{k\to\infty}a_k=a>0$. Thus $a_k\neq a$ for every $k\in\Naturals$. Additionally, there exists an $\varepsilon_a$ such that $(a-\varepsilon_a,a+\varepsilon_a)\setminus\{a\}\cap A=\emptyset$ and thus also $(a-\varepsilon_a,a+\varepsilon_a)\cap\{a_k \:|\: k\in\Naturals\}=\emptyset$, in a contradiction to the definition of a limit. Thus $a=0$, and for every $b\in A\setminus\{0\}$ the set $A\cap[b,\infty)$ must be finite. Either almost all the elements of the sequence are equal to zero, or they are all positive. In both cases, $A=\{0\}\cup\{a_k \:|\: k\in\Naturals\}$.
\end{proof}
\begin{definition}
	A totally ordered set has the order type $1+\omega^*$ if it is order isomorphic to the set $\{0\}\cup\{\tfrac{1}{2^k}\:|\: k\in\Naturals\}$ with the usual order on $\mathbb{R}$.
\end{definition}
\begin{remark}
	$\omega$ denotes the order of the natural numbers $\Naturals$. Thus, we denote by $\omega^*$, as is done in \cite{dovgoshey2022range}, the reversed order: the usual order on $\mathbb{Z}\cap(-\infty,0]$. 
\end{remark}
\begin{lemma}
	\label{orderTypeLemma}
	Let $A\subset[0,\infty)$ be a sequentially-descending-to-zero set. Then if $A$ is finite it has a finite order type, and if $A$ is infinite it has the $1+\omega^*$ order type (when joined with the usual order on $\mathbb{R}$). In both cases, every non empty subset of $A$ has a maximal element.
\end{lemma}
\begin{proof}
	Let $A\subset[0,\infty)$. The set $A$ is totally ordered by the usual order on $\mathbb{R}$.\\
	Suppose $A$ is sequentially-descending-to-zero. If $|A|\in\Naturals$, then it has necessarily a finite order type. Else, $|A|=\aleph_0$. Let $(a_k)_{k\in\Naturals}$ be a weakly monotonically decreasing sequence such that $A=\{0\}\cup\{a_k \:|\: k\in\Naturals\}$. Since $A$ is infinite, we must have $A\setminus\{0\}=\{a_k \:|\: k\in\Naturals\}$. Therefore, suppose w.l.o.g. that $(a_k)_{k\in\Naturals}$ is a strictly monotonically decreasing sequence. The function $a_k\mapsto\tfrac{1}{2^k}$ is clearly order preserving. Thus, since $0<a_k$ for every $k\in\Naturals$, the set $A$ has the order type $1+\omega^*$.\\
	Both order type satisfies the ``reverse'' well-order property: every non-empty subset of $A$ has a greatest element. 
\end{proof}
Define now the next property of metrics.
\begin{definition}[sequentially-descending-to-zero metric]
	Let $(X,d)$ be a metric space. The metric $d$ is called \emph{sequentially-descending-to-zero} if there exists a sequentially-descending-to-zero set $A\subset[0,\infty)$ such that $\image d=A$.
\end{definition} 
\begin{lemma}
	\label{openBallClosedBallLemma}
	Let $(X,d)$ be a sequentially-descending-to-zero metric space.\\
    Then every non-empty open ball in $X$ coincides with a closed ball. If $d$ is also an ultrametric, then $B(x,r)=\closedb(x,\diam B(x,r))$ for all $x\in X$ and $r\in(0,\infty)$.
\end{lemma}
\begin{proof}
	Let $x\in X$ and $r\in(0,\infty)$. Since $d$ is sequentially-descending-to-zero, by Lemma \ref{orderTypeLemma} every subset of $\image d$ has a greatest element.
	Denote
	\begin{equation*}
		r^*=\max\image d\cap[0,r)<r
	\end{equation*}
	Thus clearly $\closedb(x,r^*)\subseteq B(x,r)$.\\
	Let $y\in B(x,r)$. Then $d(x,y)\in\image d\cap[0,r)$ and thus $d(x,y)\leq r^*$, so that $y\in\closedb(x,r^*)$. Thus $B(x,r)=\closedb(x,r^*)$ and it is indeed a closed ball.\\
	Suppose now that $d$ is also ultrametric.\\
	By points \ref{diamItem} and \ref{clopenItem} of Proposition \ref{UltrametricBallsBasics} we have 
	\begin{equation*}
		\diam B(x,r)=\sup\{d(z,y)\:|\: z,y\in B(x,r)\}=\max\{d(x,y)\:|\: y\in B(x,r)\}\leq r^*,
	\end{equation*}
	and clearly $B(x,r)=\closedb(x,\diam B(x,r))$.
\end{proof}
Similarly,
\begin{lemma}
	\label{closedBallopenBallLemma}
	Let $(X,d)$ be a sequentially-descending-to-zero metric space. Then for every point $x\in X$ and every radius $r\in(0,\infty)$ there exists a radius $\rho_{x,r}\in(r,\infty)$ such that
	\begin{equation*}
		\closedb(x,r)=B(x,\rho_{x,r}).
	\end{equation*}
\end{lemma}
\begin{proof}
	Let $x\in X$ and $r\in(0,\infty)$. Since $d$ is sequentially-descending-to-zero, the set $\image d$ is by definition sequentially-descending-to-zero. Thus by Proposition \ref{SDTZsetProperty}, there exists an $\varepsilon_r\in(0,\infty)$ such that
    \begin{equation*}
        (r-\varepsilon_r,r+\varepsilon_r)\setminus\{r\}\cap\image d=\emptyset.
    \end{equation*}
    Denote $\rho_{x,r}=r+\tfrac{\varepsilon_r}{2}$. Thus
    \begin{equation*}
		\closedb(x,r)=B\left(x,r+\tfrac{\varepsilon_r}{2}\right)=\closedb(x,\rho_{x,r}).
    \end{equation*}
\end{proof}
The proof of Proposition \ref{universalProposition} is based on the one from \cite[Lemma 4.3]{dovgoshey2022range}. The importance of this proposition comes from the topological consequences that can be learned about ultrametrizable spaces when analyzing sequentially-descending-to-zero ultrametrics.
\begin{proposition}
	\label{universalProposition}
	Let $(X,\mathcal{O})$ be an ultrametrizable topological space. There exists a sequentially-descending-to-zero ultrametric that induces the topology $\mathcal{O}$.
\end{proposition}
\begin{proof}
	Let $(X,\mathcal{O})$ be an ultrametrizable topological space and let $d$ an ultrametric that induces this topology. Denote the weakly monotonically increasing function
	\begin{equation*}
		g\colon[0,\infty)\longrightarrow\overline{\left\{\tfrac{1}{2^k}\:|\: k\in\Naturals\right\}}
	\end{equation*}
	by
	\begin{equation}
		g(x)=
		\begin{cases}
			1, & \text{if }x\geq 1,\\
			2^{\lfloor\log_2x\rfloor}, & \text{if }x\in(0,1),\\
			0, & \text{if }x=0.\\
		\end{cases}
	\end{equation}
	such that if $\tfrac{1}{2^k}\leq x<\tfrac{1}{2^{k-1}}$ for some $k\in\Naturals\setminus\{0\}$ then $g(x)=\tfrac{1}{2^k}$. Define the ultrametric
	\begin{equation*}
		\widetilde{d}\colon X\times X\longrightarrow\overline{\left\{\tfrac{1}{2^k}\:|\: k\in\Naturals\right\}}
	\end{equation*}
	by
	\begin{equation}
		\widetilde{d}=g\circ d,
	\end{equation}
	which is indeed an ultrametric: $g^{-1}(0)=0$ so that $\widetilde{d}(x,y)=0\:\Leftrightarrow\:x=y$; $\widetilde{d}$ is clearly symmetric because $d$ is; and if $x,y,z\in X$ then since $g$ is weakly monotonically increasing
	\begin{equation*}
		\begin{split}
			\widetilde{d}(x,y)&=g(d(x,y))\leq g\left(\max\{d(x,z),d(y,z)\}\right)=\\
			& =\max\{g(d(x,z)),g(d(y,z))\}=\max\{\widetilde{d}(x,z),d(y,z)\}.\\
		\end{split}
	\end{equation*}
	Additionally, $\widetilde{d}$ also induces the same topology on $X$:\\
	Let $(x_k)_{k\in\Naturals}\in X^\Naturals$ such that $\lim_{k\to\infty}x_k=x\in X$ in the ultrametric space $(X,d)$. Then $\lim_{k\to\infty}d(x_k,x)=0$ and thus $\lim_{k\to\infty}\widetilde{d}(x_k,x)=0$ and hence $\lim_{k\to\infty}x_k=x\in X$ also in the ultrametric space $(X,\widetilde{d})$.\\
	Finally, the ultrametric $\widetilde{d}$ is also sequentially-descending-to-zero, as
	\begin{equation*}
		\image\widetilde{d}\subseteq\image g=\overline{\left\{\tfrac{1}{2^k}\:|\: k\in\Naturals\right\}}
	\end{equation*}
\end{proof}
\begin{comment}
    When an ultrametric is sequentially-descending-to-zero, this paper proposes to call it simply a \emph{dendronic} ultrametric, for reasons detailed in Section \ref{TreesSection}.
\end{comment}
\begin{comment}
	The dendronic ultrametrics have an important role in characterizing compact ultrametrizable space, as was shown in \cite[Theorem 4.7]{dovgoshey2022range}:
	\begin{proposition}
		Let $(X,\mathcal{O})$ be an ultrametrizable space. Then $(X,\mathcal{O})$ is compact if and only if every ultrametric $d$ compatible with the topology $\mathcal{O}$ is dendronic.
	\end{proposition}
	\begin{proof}
		Let $(X,d)$ be a compact ultrametric space. Suppose by contradiction that $d$ is not dendronic. Thus, $\image d$ is not sequentially-descending-to-zero. Thus there exists an $r\in(0,\infty)$ such that $(r-\varepsilon,r+\varepsilon)\setminus\{r\}\cap\image d\neq\emptyset$ for every $\varepsilon>0$. Suppose w.l.o.g. that $(r,r+\varepsilon)\cap\image d\neq\emptyset$ for every $\varepsilon>0$ (the proof for the other case is analogous). Suppose that $r\in\image d$, and let $x,y\in X$ such that $d(x,y)=r$. By point \ref{CompactClopen} of Proposition \ref{UltrametricBallsBasics}, the closed ball $\closedb(x,r)$ is also an open ball. Since $\diam\closedb(x,r)$, there must exist a $q\in(r,\infty)$ such that $B(x,q)=\closedb(x,r)$. 
	\end{proof}
\end{comment}

\subsection{Vitali's covering lemma}
The classical infinite version of Vitali's covering lemma states that given an arbitrary set of open balls with uniformly bounded radii in a metric space, there exists a subset of mutually disjoint balls (out of the original set), such that when increasing their radii times 4, they cover all the original set of balls (see, for example, \cite[Theorem 8.1]{mauldin2009graph}). Additionally, when the space is also separable, there exists a countable subset that satisfies this requirements. In the case of a sequentially-descending-to-zero ultrametric, Vitali's covering lemma is satisfied in an even stronger version: there is no need to increase the radii of any ball to cover the original set.
\begin{lemma}[Vitali's covering lemma]
    \label{VitaliLemma}
    Let $(X,d)$ be a sequentially-descending-to-zero ultrametric space, and let $\mathcal{C}$ be a set of open balls. There exists a subset $\mathcal{C}'\subseteq\mathcal{C}$ of mutually disjoint balls such that
    \begin{equation}
        \bigcup_{B\in\mathcal{C}}B=\biguplus_{B\in\mathcal{C}'}B.
    \end{equation}
    If the ultrametric space is also separable, or alternatively if there exists a Borel probability measure $\mu$ such that $\mu(B)>0$ for every $B\in\mathcal{C}'$, then the subset of balls $\mathcal{C}'$ is countable.
\end{lemma}
\begin{proof}
    Let $(X,d)$ be a sequentially-descending-to-zero ultrametric space and let $\mathcal{C}$ be a set of open balls. Note that $\max\{\diam B \:|\: B\in\mathcal{C}\}\leq\max\image d$.\\
    By lemma \ref{openBallClosedBallLemma}, every open ball in $\mathcal(C)$ is also closed ball.
    Let $(a_k)_{k\in\Naturals}$ be a weakly monotonically decreasing sequence that converges to zero and such that $\image d = \{0\}\cup\{a_k \:|\: k\in\Naturals\}$. Suppose w.l.o.g that if $a_k>0$ then $a_k>a_{k+1}$ for every $k\in\Naturals$.
    Denote
    \begin{equation*}
        M=\sup\{\diam B\:|\: B\in\mathcal{C}\}\in\image d.
    \end{equation*}
    There is a natural number $K_0$ such that $a_{K_0}=M$. Define
    \begin{equation*}
        \mathcal{C}_0=\{B\in\mathcal{C}\:|\:\diam B=M=a_{K_0}\}.
    \end{equation*}
    Since all the balls in $\mathcal{C}$ are closed ball with the same radius, they all must be disjoint by point \ref{egocentricity} of Proposition \ref{UltrametricBallsBasics}. By point \ref{inclusionBalls} of Proposition \ref{UltrametricBallsBasics}, we get that if another ball of a smaller radius intersects one of the balls in $\mathcal{C}_0$, then it is included in it. For every $k\in\Naturals\setminus\{0\}$ define
    \begin{equation*}
        \mathcal{C}_k=\left\{B\in\mathcal{C}\:\middle|\:\diam B=a_{K_0+k}\text{ and }B\cap\biguplus_{B'\in\biguplus^{k-1}_{j=0}\mathcal{C}_j}B'=\emptyset\right\},
    \end{equation*}
    and finally denote
    \begin{equation*}
        \mathcal{C}'=\biguplus_{k\in\Naturals}\mathcal{C}_k,
    \end{equation*}
    which is by its construction a set of mutually disjoint balls. We also have
    \begin{equation*}
        \bigcup_{B\in\mathcal{C}}B=\biguplus_{B\in\mathcal{C}'}B,
    \end{equation*}
    since if there exists a ball $B\in\mathcal{C}$, then
    \begin{equation*}
        \closedb\left(x,\max\{a_{K_0+k}\:|\: k\in\Naturals\text{ and }\closedb(x,a_{K_0+k})\in\mathcal{C}\}\right)\in\mathcal{C}'
    \end{equation*}
    for every $x\in B$.\\
    Suppose now that $(X,d)$ is also separable. Let $A\subseteq X$ be a countable dense subset. Let $B\in\mathcal{C}'$. Since $B$ is also an open ball, we must have $B\cap A\neq\emptyset$. Moreover, $B=\closedb(x,\diam B)$, where $x\in B\cap A$ is arbitrary. Thus we must have $|\mathcal{C}'|\leq|A|\leq\aleph_0$.\\
    Suppose now alternatively that $(X,d)$ is not necessarily separable, but that there exists a Borel probability measure $\mu$ such that $\mu(B)>0$ for every $B\in\mathcal{C}'$. Suppose by contradiction that $|\mathcal{C}'|>\aleph_0$. Then, by the axiom of countable choice, there must exist some $k\in\Naturals\setminus\{0\}$ such that the set
    \begin{equation*}
        \mathcal{C}'_k=\left\{B\in\mathcal{C}'\:\middle|\:\tfrac{1}{k+1}<\mu(B)\leq\tfrac{1}{k}\right\}
    \end{equation*}
    is uncountable. Denote by $\mathcal{C}''_k\subset\mathcal{C}'_k$ a strict subset such that $|\mathcal{C}''_k|=\aleph_0$. Then
    \begin{equation*}
        1\geq\mu\left(\biguplus_{B\in\mathcal{C}'}B\right)\geq\mu\left(\biguplus_{B\in\mathcal{C}''_k}B\right)\stackrel{|\mathcal{C}''_k|=\aleph_0}{=}\sum_{B\in\mathcal{C}''_k}\mu(B)\geq\sum_{i=0}^\infty\tfrac{1}{k+1}=\infty,
    \end{equation*}
    a contradiction. Thus $|\mathcal{C}'|\leq\aleph_0$.
\end{proof}

\section{Pruned trees and branches over alphabets}
\label{Trees}
The relationship between trees and ultrametric spaces has been a subject of study for some time. While some researchers have explored the connections between ultrametric spaces and $\mathbb{R}$-trees (e.g., \cite{hughes2004trees}), this paper adopts the discrete tree approach similar to that employed in \cite{brian2015completely}.

This section demonstrates that sequentially-descending-to-zero ultrametric spaces are isometric to spaces of branches of pruned trees. Leveraging Proposition \ref{universalProposition}, proven in \cite[Lemma 4.3]{dovgoshey2022range}, we establish that every ultrametrizable space is homeomorphic to a space of branches (referred to in this paper as clade spaces) of a pruned tree. This result generalizes the well-known theorem concerning completely ultrametrizable spaces (see \cite[Proposition 2.1]{brian2015completely}).

Importantly, topological properties such as completeness, separability, and compactness manifest in straightforward ways within tree structures. This simplification facilitates the identification and analysis of these properties in ultrametrizable spaces, offering a powerful framework for their study.

\begin{definition}[tree]
	Let $\mathcal{A}$ be a non-empty set,  called an \emph{alphabet}. A tree over the alphabet $\mathcal{A}$ is a non-empty subset $T\subseteq\mathcal{A}^{<\omega}=\bigcup_{k\in\Naturals}\mathcal{A}^k$ that satisfies 
	\begin{enumerate}
		\label{TreeAxioms}
		\item[T.1)] If $p\in T$ and $p'$ is a prefix of $p$, then also $p'\in T$.
	\end{enumerate}
	If $T$ also satisfies the condition
	\begin{enumerate}
		\item[T.2)] For every $p\in T$ there exists a $p'\in T$ such that $p$ is a strict prefix of $p'$,
	\end{enumerate}
	then it is said that it is \emph{pruned}. An element of $T$ is called a \emph{position}.\\
	The strict partial order $\prec\subset T\times T$ is defined by $p\prec p'$ when $p$ is a strict prefix of $p'$, and in that case $p$ is said to be an \emph{ascendant} of $p'$ and the latter a \emph{descendant} of the former. An immediate ascendant (prefix sequence with length smaller by one) is called a \emph{parent} and an immediate descendant a \emph{child} or an \emph{offspring}.
\end{definition}
The elements of a tree, \emph{i.e.}, its positions, are thus finite sequences of elements from the adequate alphabet. The tree has a \emph{root}, the empty set/sequence, which is a prefix of all positions of the tree. Note that this definition of a tree is a descriptive set theoretic one. The length of a position $p$ (\emph{i.e.}, the length of the finite sequence $p$) is denoted by $\len(p)$.
\begin{definition}[branches, canopy]
	\label{branchCanopyDefinition}
	Let $T$ be a tree over the alphabet $\mathcal{A}$. A branch of $T$ is an infinite sequence in $\mathcal{A}^\Naturals$ such that all its finite prefixes are in $T$.\par
	The set of all branches of the tree $T$ is denoted by $\llbracket T\rrbracket$ and is called the \emph{end space}, \emph{branch space}, \emph{body}, or simply \emph{canopy} of the tree $T$.
\end{definition}
\begin{remark}
	A branch is a maximal linearly ordered subset of $T$. 
\end{remark}
\begin{example}[complete binary tree]
	The \emph{complete binary tree} is the pruned tree over the alphabet $\mathcal{A}=\{0,1\}$ defined by
	\begin{equation*}
		T=\mathcal{A}^{<\omega}=\bigcup_{k\in\Naturals}\mathcal{A}^k,
	\end{equation*}
	and its canopy is the set of all binary sequences, \emph{i.e.}, $\llbracket T\rrbracket=\Dya$.
\end{example}
Two additional important definitions in the realm of trees are the followings:
\begin{definition}[the subtree $T_p$]
	\label{subtreeDefinition}
	Let $T$ be a tree over an alphabet $\mathcal{A}$, and let $p\in T$. The subtree $T_p\subseteq T$ is defined by
	\begin{equation*}
		T_p=\{p'\in T \:|\: p\preccurlyeq p'\text{ or }p'\preccurlyeq p\}.
	\end{equation*}
\end{definition}
\begin{definition}[tree topology of a canopy]
\label{CladeCanopyDefinition}
	Let $T$ be a pruned tree over the alphabet $\mathcal{A}$. The \emph{tree topology} over the canopy $\llbracket T\rrbracket$, denoted by $\widetilde{\mathcal{O}}_T$, is the topology generated by the basis
	\begin{equation*}
		\left\{\llbracket T_p\rrbracket\:|\: p\in T\right\}.
	\end{equation*}
	If $\emptyset\neq P\subseteq\llbracket T\rrbracket$ is a subset of branches, then the tree topology over $P$ is the subspace topology induced by $\widetilde{\mathcal{O}}_T$ and is denoted by $\widetilde{\mathcal{O}}_{T,P}$. We call $(P,\widetilde{\mathcal{O}}_{T,P}) $ a \emph{clade space}\footnote{From $\kappa\lambda\alpha\delta o\sigma$ (``Klados'') which means ``branch'' in ancient Greek.} of the pruned tree $T$.
\end{definition}
\begin{remark}
	The tree topology over the canopy $\llbracket T\rrbracket\subseteq\mathcal{A}^\mathbb{N}$ is the same as the subspace topology inherited from the space $\mathcal{A}^\mathbb{N}$ endowed with the infinite product of the discrete topology over $\mathcal{A}$.
\end{remark}
\begin{notation}
	\label{BinaryTreeSpecialNotation}
	We will sometimes denote $\llbracket p\rrbracket$ instead of $\llbracket T_p\rrbracket$ for positions $p\in T$.
\end{notation}
From this point on, this paper will consider only pruned trees and will omit sometimes  the adjective ``pruned'' to make the writing lighter. Similarly, ``subtree'' would mean a pruned subtree, unless explicitly mentioned otherwise.
\begin{definition}
	\label{SubcanopyDefinition}
	Let $T$ be a pruned tree over the alphabet $\mathcal{A}$ and let $T'\subseteq T$ be a pruned subtree. Then the canopy $\llbracket T'\rrbracket$ is said to a \emph{subcanopy} of $\llbracket T\rrbracket$.
\end{definition}
\begin{proposition}[tree topology is ultrametrizable]
	\label{TreeTopologyUltrametrizable}
	Let $T$ be a pruned tree over the alphabet $\mathcal{A}$ and let $\emptyset\neq P\subseteq\llbracket T\rrbracket$. There exists a sequentially-descending-to-zero ultrametric $d$ such that $d$ induces the tree topology $\widetilde{\mathcal{O}}_{T,P}$. The ultrametric space $(P,d)$ is complete if and only if $P$ is a subcanopy of $\llbracket T\rrbracket$ (see Definition \ref{SubcanopyDefinition}).
\end{proposition}
\begin{proof}
	The dyadic ultrametric over $\mathcal{A}^\Naturals$ is defined by
	\begin{equation}
		\label{dyadicDendronicUltrametric}
		d_2(x,y)=
		\begin{cases}
			0, & \text{if}\quad x=y,\\
			2^{-\min\{k\in\Naturals\:|\:x_k\neq y_k\}}, & \text{else}.\\
		\end{cases}
	\end{equation}
	We show it is indeed an ultrametric:\\
	Let $x,y,z\in\mathcal{A}^\Naturals$. Denote
	\begin{equation*}
		\begin{split}
			& k_{x,y} = \min\{k\in\Naturals\:|\:x_k\neq y_k\}\\
			& k_{y,z} = \min\{k\in\Naturals\:|\:y_k\neq z_k\}\\
			& k_{x,z} = \min\{k\in\Naturals\:|\:x_k\neq z_k\}
		\end{split}
	\end{equation*}
	If $k_{x,z}\geq k_{x,y}$, then for every $0\leq i<k_{x,y}$ we have both  $x_i=z_i$ and $x_i=y_i$ so by transitivity also $y_i=z_i$ and thus $k_{y,z}\geq k_{x,y}$. But since $z_{k_{x,y}}=x_{k_{x,y}}\neq y_{k_{x,y}}$ we have $k_{y,z}=k_{x,y}$. Else, if $k_{x,z}<k_{x,y}$, then $z_{k_{x,z}}\neq x_{k_{x,z}}=y_{k_{x,z}}$ and thus $k_{y,z}\leq k_{x,z}<k_{x,y}$. In both cases we have
	\begin{equation*}
		\begin{split}
			k_{x,y}&\geq\min\{k_{x,z},k_{y,z}\}\\
			\Rightarrow\:\:-k_{x,y}&\leq -\min\{k_{x,z},k_{y,z}\}=\max\{-k_{x,z},-k_{y,z}\}\\
			\Rightarrow\:\:2^{-k_{x,y}}&\leq 2^{\max\{-k_{x,z},-k_{y,z}\}}=\max\{2^{-k_{x,z}},2^{-k_{y,z}}\}\\
			\Rightarrow\:\: d_2(x,y)&\leq\max\{d_2(x,z),d_2(y,z)\}
		\end{split}
	\end{equation*}		
	and since $\image d_2=\overline{\left\{\tfrac{1}{2^k}\:|\: k\in\Naturals\right\}}$ it is clearly a sequentially-descending-to-zero ultrametric on $\mathcal{A}^\Naturals$. Since $\llbracket T\rrbracket\subseteq\mathcal{A}^\Naturals$ and $P\subseteq\llbracket T\rrbracket$, $d_2$ is also a sequentially-descending-to-zero ultrametric over $P$ (when restricted to $P\times P$). This ultrametric induces the tree topology since
	\begin{equation*}
		\llbracket T_p\rrbracket = B\left(x_p,\tfrac{1}{2^{\len(p)}}\right)
	\end{equation*}
	for every position $p\in T$, where $x_p\in\llbracket T\rrbracket$ is a branch such that $p$ is a finite prefix of it.\\
	Suppose $P=\llbracket T'\rrbracket$ where $T'\subseteq T$ is a subtree. Then ($\llbracket T'\rrbracket,d_2)$ is complete:\\
	Let $(x^k)$ be a Cauchy sequence of branches in the subcanopy $\llbracket T'\rrbracket$. For every $l\in\Naturals$, let $k_l\in\Naturals$ such that $d_2(x^i,x^j)<\tfrac{1}{2^l}$ for every $i,j\geq k_l$. We now construct a sequence of elements in $T'$ (finite sequences) with increasing lengths: For every $l\in\Naturals$ denote $y^l=(x^{k_l}_m)_{m=0}^l\in T'$.\\
    The sequence $y=(y^l)_{l\in\Naturals}$ is a sequence of position such that each one is the offspring of the precedent position, and thus it is a branch of the subtree $T'$. Moreover, for every $l\in\Naturals$ the elements of $y^l$ equal the $l+1$ first elements of the sequence $x^i$ for every $i\geq k_l$. The branch $y\in\llbracket T'\rrbracket$ is the limit of $(x^k)_{k\in\Naturals}$: Let $\varepsilon>0$. Denote
	\begin{equation*}
		K_\varepsilon=K_{\max\left\{0,\left\lfloor\log_2\left(\tfrac{1}{\varepsilon}\right)\right\rfloor+1\right\}}.
	\end{equation*}
	Then for every $k\geq K_\varepsilon$ we have
	\begin{equation*}
		d_2(y,x^k)<\frac{1}{2^{\left\lfloor\log_2\left(\tfrac{1}{\varepsilon}\right)\right\rfloor+1}}<\varepsilon.
	\end{equation*}
	Suppose now that $(P,d_2)$ is complete. We show there exists a subtree $T'\subseteq T$ such that $P=\llbracket T'\rrbracket$:\\
	We define
	\begin{equation*}
		T'=\{\emptyset\}\cup\{(x_i)_{j=0}^k \:|\: k\in\Naturals\text{ and }(x_i)_{i\in\Naturals}\in P\}
	\end{equation*}
	clearly $P\subseteq\llbracket T'\rrbracket$. Let $x\in\llbracket T'\rrbracket$. Then there exists a sequence of positions in $T'$ with increasing lengths $(p^k)_{k\in\Naturals}$ such that $p^k$ is a prefix of $x$ for every $k\in\Naturals$. For every $k\in\Naturals$, let $x^k\in P$ be an infinite sequence such that $p^k$ is a finite prefix $x^k$ (at least one such infinite sequence exists by the definition of $T'$, note that the axiom of countable choice is employed here). Note that
	\begin{equation*}
		d_2(x^k,x)\leq\frac{1}{2^{\len(p^k)}}
	\end{equation*}
	for every $k\in\Naturals$. Thus $\lim_{k\to\infty}x^k=x$ in $\mathcal{A}^\Naturals$. But by the strong triangle inequality we also get that $(x^k)_{k\in\Naturals}$ is Cauchy. Thus, by the assumption that $(P,d_2)$ is complete we get that $x\in P$. Therefore $\llbracket T'\rrbracket\subseteq P$, and we get that $P=\llbracket T'\rrbracket$.
\end{proof}
\begin{comment}
    \begin{remark}
    	When $\Dya$ is considered as the canopy of the complete binary tree, the dyadic ultrametric defined in Section \ref{dyadicMetricSection} is the same as defined in the proof of Proposition \ref{TreeTopologyUltrametrizable}. Therefore, the dyadic ultrametric space of binary sequences is complete, as a special case of Proposition \ref{TreeTopologyUltrametrizable}.
    \end{remark}
\end{comment}
\begin{corollary}
	\label{ClosedSubcanopyCorollary}
	Let $T$ be a pruned tree over the alphabet $\mathcal{A}$ and let $\emptyset\neq P\subseteq\llbracket T\rrbracket$. Then $P$ is closed (according to the tree topology) if and only if it is a subcanopy.
\end{corollary}
\begin{proof}
	By Proposition \ref{TreeTopologyUltrametrizable}, the canopy $\llbracket T\rrbracket$ with the tree topology is completely metrizable. Since a subset of a completely metric space is closed if and only if it is complete, it is sufficient to show that $\emptyset\neq P\subseteq\llbracket T\rrbracket$ is complete if and only if it is a subcanopy, which was shown in Proposition \ref{TreeTopologyUltrametrizable}.
\end{proof}
The converse of Proposition \ref{TreeTopologyUltrametrizable} is also true: Corollary \ref{HomeomorphismCorollary} claims that every ultrametrizable space is homeomorphic to a clade space of a pruned tree. 
\begin{definition}[representing clade spaces]
	\label{representingDefinition}
	Let $(X,d)$ be a sequentially-descending-to-zero ultrametric space. Let $(r_k)_{k\in\Naturals}$ be a weakly monotonically decreasing sequence of non-negative real numbers such that w.l.o.g.\footnote{Given a weakly monotonically decreasing sequence of non-negative real numbers $(a_k)_{k\in\Naturals}$, define a sequence $(b_k)_{k\in\Naturals}$ by $b_k=\max\left(\{a_j \:|\ :j>k\}\setminus\{a_k\}\right)$ if $a_k>0$ and else $b_k=0=a_k$.} for every $k\in\Naturals$, if $r_k>0$ then $r_k>r_{k+1}$, and also such that
	\begin{equation*}
		\image d = \{0\}\cup\{r_k \:|\: k\in\Naturals \}.
	\end{equation*}
	Thus, $(r_k)_{k\in\Naturals}$ is either a strictly monotonically decreasing sequence of positive real numbers that tends to 0, or alternatively a sequence that is a concatenation of a strictly monotonically decreasing finite sequence of positive numbers with a constant infinite sequence of 0s. Define the tree $T_{(X,d)}$ over the alphabet
	\begin{equation}
		\mathcal{A}_{(X,d)}=\left\{\closedb(x,r_k)\:|\: x\in X\text{ and }k\in\Naturals\setminus\{0\}\right\}
	\end{equation}
	by
	\begin{equation}
		T_{(X,d)}=\{\emptyset\}\cup\left\{\left(\closedb(x,r_j)\right)_{j=1}^k \:\middle|\: x\in X\text{ and }k\in\Naturals\setminus\{0\}\right\},
	\end{equation}
	the metric $d_{T_{(X,d)}}\colon\mathcal{A}^\Naturals\times\mathcal{A}^\Naturals$$\longrightarrow\image d$ by
	\begin{equation*}
		d_{T_{(X,d)}}(x,y)=
		\begin{cases}
			0, & \text{if } x=y,\\
			r_{\min\{k\in\Naturals\:|\:x_k\neq y_k\}}, & \text{else}.\\
		\end{cases}
	\end{equation*}
	and the set of branches
	\begin{equation}
		\llceil X\rrceil = \left\{\left(\closedb(x,r_k)\right)_{k\in\Naturals}\:\middle|\: x\in X\right\}.
	\end{equation}
	We say that the tree $T_{(X,d)}$ and the clade space $\left(\llceil X\rrceil,d_{T_{(X,d)}}\right)$ represent the ultrametric space $(X,d)$.
\end{definition}
Note that the same tree (but not the same clade space) can represent many different ultrametric spaces. Additionally, since $d$ is an ultrametric, $\closedb(x,r_k)=\closedb(y,r_k)\in\mathcal{A}_{(X,d)}$ for every $k\in\Naturals$ and every $y,x\in X$ with $d(x,y)\leq r_k$. Note that the root of the tree, \emph{i.e.}, the empty set, represents the whole space $X$, and that the positions of length $k$ represent the closed balls with radius $r_k$ for every positive natural number $k$.
\begin{theorem}
	\label{IsometricToTreeTheorem}
	Let $(X,d)$ be a sequentially-descending-to-zero metric space. Then $(X,d)$ is isometric to its representing clade space $\left(\llceil X\rrceil,d_{T_{(X,d)}}\right)$.
\end{theorem}
\begin{proof}
	Define the function $\phi\colon X\longrightarrow\llceil X\rrceil$ by
	\begin{equation}
		\phi(x)=\left(\closedb(x,r_k)\right)_{k\in\Naturals},
	\end{equation}
	whose inverse is equal to
	\begin{equation*}
		\phi^{-1}((\closedb_k)_{k\in\Naturals})=\bigcap_{k\in\Naturals}\closedb_k,
	\end{equation*}
	when identifying a singleton with its element.\\
	Let $x,y\in X$ and let $k_{x,y}\in\Naturals$ such that $d(x,y)=r_{k_{x,y}}$. Thus for every $0\leq j \leq k_{x,y}$ we have $\closedb(x,r_j)=\closedb(y,r_j)$ and thus
	\begin{equation*}
		\begin{split}
			&d_{T_{(X,d)}}(\phi(x),\phi(y))=d_{T_{(X,d)}}((\closedb(x,r_j))_{j\in\Naturals},(\closedb(y,r_j))_{j\in\Naturals})=\\
			&=r_{\min\{j\in\Naturals\:|\: \closedb(x,r_{j+1})\neq \closedb(y,r_{j+1})\}}=r_{k_{x,y}}=d(x,y).
		\end{split}
	\end{equation*}
	Let $(B_i)_{i\in\Naturals},(D_i)_{i\in\Naturals}\in\llceil X\rrceil$ such that $d_{T_{(X,d)}}((B_i)_{i\in\Naturals},(D_i)_{i\in\Naturals}) = r_{k^*}$ for some $k^*\in\Naturals$. Thus $B_i\neq D_i$ for every $k^*\leq i$. Denote $x=\phi^{-1}((B_i)_{i\in\Naturals})$ and $y=\phi^{-1}((D_i)_{i\in\Naturals})$. Then
	\begin{equation*}
		\closedb(x,r_{i+1})=B_i\text{  and  }\closedb(y,r_{i+1})=D_i
	\end{equation*}
	for every $i\in\Naturals$. Hence
	\begin{equation*}
		\closedb(x,r_{k^*+1})\neq\closedb(y,r_{k^*+1})
	\end{equation*}
	but
	\begin{equation*}
		\closedb(x,r_{k^*})=\closedb(y,r_{k^*})
	\end{equation*}
	and therefore \begin{equation*}
		d(\phi^{-1}((B_i)_{i\in\Naturals}),\phi^{-1}((D_i)_{i\in\Naturals}))=d(x,y)=r_{k^*}=d_{T_{(X,d)}}((B_i)_{i\in\Naturals},(D_i)_{i\in\Naturals})
	\end{equation*}
	and hence $\phi$ is indeed an isometry.
\end{proof}
\begin{corollary}
	\label{HomeomorphismCorollary}
	An ultrametrizable space is homeomorphic to a clade space of a pruned tree.
\end{corollary}
\begin{proof}
	Let $(X,\mathcal{O})$ be an ultrametrizable topological space.\\
	Let $d$ be a sequentially-descending-to-zero ultrametric over $X$ such that $\mathcal{O}=\mathcal{O}_{(X,d)}$, which exists by Proposition \ref{universalProposition}. By Theorem \ref{IsometricToTreeTheorem}, $(X,d)$ is isometric to the representing clade space $(\llceil X\rrceil,d_{T_{(X,d)}})$, and thus $(X,\mathcal{O})$ is homeomorphic to a clade space.
\end{proof}
Corollary \ref{HomeomorphismCorollary} is useful because some topological properties have very simple characterizations in the case of clade spaces. The rest of this section details some examples.

\subsection{Completeness}
As proved in Proposition \ref{TreeTopologyUltrametrizable}, a clade space of a pruned tree is complete if and only if it is a subcanopy. Consequently, completeness has a simple characterization using the representing clade space of an ultrametrizable space.
\begin{proposition}
    Let $(X,\mathcal{O})$ be an ultrametrizable space.\\
    Let $T_X$ and $(\llceil X\rrceil,d_{T_X})$ be a representing tree and a corresponding representing clade space of $(X,\mathcal{O})$. Then the completion of the clade space $(\llceil X\rrceil,d_{T_X})$ is the canopy $(\llbracket T_X\rrbracket,d_{T_X})$  (and therefore $(X,\mathcal{O})$ is completely ultrametrizable if and only if $\llceil X\rrceil=\llbracket T_X\rrbracket$).
\end{proposition}
\begin{proof}
	Let $(X,\mathcal{O})$ be an ultrametrizable space. Let $T_X$ and $(\llceil X\rrceil,d_{T_X})$ be a representing tree and a corresponding representing clade space of $(X,\mathcal{O})$, which exist by Theorem \ref{IsometricToTreeTheorem}. Let $(r_k)_{k\in\Naturals}$ be a weakly monotonically decreasing sequence such that if $r_k>0$ then $r_{k+1}<r_k$ for every $k\in\Naturals$ and also 
	\begin{equation*}
		\image d_{T_X} = \{0\}\cup\{r_k \:|\: k\in\Naturals \}.
	\end{equation*}
	By Proposition \ref{TreeTopologyUltrametrizable}, the ultrametric space $\left(\llbracket T_X\rrbracket,d_{T_X}\right)$ is a complete ultrametric space. Thus if $\llceil X\rrceil=\llbracket T_X\rrbracket$ then $(X,\mathcal{O})$ is completely ultrametrizable.\\
	Suppose now that $(X,\mathcal{O})$ is completely ultrametrizable. Then by Proposition \ref{TreeTopologyUltrametrizable}, $\llceil X\rrceil$ must be a subcanopy, \emph{i.e.}, there exists a non-empty subtree $T'\subseteq T_X$ such that $\llceil X\rrceil=\llbracket T'\rrbracket$. Suppose by contradiction that $T'\subset T_X$. Let $y\in\llbracket T_X\rrbracket\setminus\llbracket T'\rrbracket$. Thus there is a sequence of closed balls $(y_i)_{i\in\Naturals}$ such that there is no element $x\in X$ such that $y_i=\closedb(x,r_{i+1})$ for every $i\in\Naturals$. Let $(x^i)_{i\in\Naturals}$ be a sequence of elements in $X$ such that $y_i=\closedb(x^i,r_{i+1})$ for every $i\in\Naturals$, which must exist by definition of the representing tree. Then $(x^i)_{i\in\Naturals}$ is a Cauchy sequence with no limit in $X$, a contradiction. Thus $\llceil X\rrceil = \llbracket T_X\rrbracket$. Consequently, there is no strict subcanopy of $\llbracket T\rrbracket$ that includes $\llceil X\rrbracket$, which makes $\llbracket T\rrbracket$ the completion of $\llceil X\rrceil$.
\end{proof}

\subsection{Separability}
\label{SeparabilityCanopySubsection}
Since separability of a metric space is equivalent to being Lindel\"of or second-countable is has a simple characterization in the case of clade spaces.
\begin{proposition}
	\label{SeparabilityCanopyProposition}
	A canopy of a pruned tree is separable (according to the tree topology) if and only if the tree has countably many positions. 
\end{proposition}
\begin{proof}
    Let $T$ be a pruned tree. Suppose that $|T|=\aleph_0$.\\
    Thus the base $\left\{\llbracket T_p\rrbracket\:|\: p\in T\right\}$, which generates the tree topology, is countable. Hence The canopy endowed with the tree topology, \emph{i.e.}, $\left(\llbracket T\rrbracket,\widetilde{\mathcal{O}}_T\right)$,  is second-countable and thus separable.\\
	Suppose that $|T|>\aleph_0$. Then, by the axiom of countable choice, there must exist at least one $k\in\Naturals\setminus\{0\}$ such that
	\begin{equation*}
		|\{p\in T \:|\: \len(p)=k\}|>\aleph_0
	\end{equation*}
	and therefore the open cover
	\begin{equation*}
		\left\{\llbracket T_p\rrbracket\:|\: \len(p)=k\right\}
	\end{equation*}
	does not have a countable subcover. Hence $(\llbracket T\rrbracket,\widetilde{\mathcal{O}}_T)$ is not Lindel\"of and thus not separable.
\end{proof}
\begin{remark}
	By Proposition \ref{SeparabilityCanopyProposition}, since the complete binary tree is countable, its canopy is separable. Thus the dyadic ultrametric space of binary sequences $(\Dya,d_2)$ is separable, a fact that can be deduced also since it is a doubling metric space (see Proposition \ref{DoublingImpliesBoundedOffsprings} below).
\end{remark}
\begin{corollary}
	\label{separableUltrametrizableCorollary}
	An ultrametrizable space is separable if and only if every pruned tree that represents it is countable.
\end{corollary}
\begin{proof}
	Let $(X,\mathcal{O})$ be an ultrametrizable space and let $T_X$ a tree that represents it. If $|T_X|=\aleph_0$, then $(\llbracket T_X\rrbracket,\widetilde{\mathcal{O}}_{T_X})$ is separable by Proposition \ref{SeparabilityCanopyProposition}, and thus also the clade space representing $(X,\mathcal{O})$ and thus also $(X,\mathcal{O})$.\\
	Suppose that $(X,\mathcal{O})$ is separable. Then also $(\llceil X\rrceil,\widetilde{\mathcal{O}}_{T_X})$ is separable and thus also $(\llbracket T_X\rrbracket,\widetilde{\mathcal{O}}_{T_X})$ since a completion of a separable metric space is separable. Hence by Proposition \ref{SeparabilityCanopyProposition}, $|T_X|=\aleph_0$.
\end{proof}

\subsection{Total boundedness and compactness}
The connection between balls and subcanopies enables a simple characterization of total boundedness of ultrametrizable spaces.
\begin{proposition}
	\label{CompactTree}
	Let $T$ be a pruned tree. Then its canopy endowed with the tree topology is compact if and only if every position of the tree has a finite number of offsprings.
\end{proposition}
\begin{proof}
	Let $T$ be a pruned tree and let $d_2$ be the dyadic metric over the canopy $\llbracket T\rrbracket$. Suppose that every position of $T$ has a finite number of offsprings. Let $\varepsilon>0$ and denote $K_\varepsilon=\max\left\{0,\left\lfloor\tfrac{1}{\varepsilon}\right\rfloor+1\right\}$. Then the collection of open balls
	\begin{equation*}
		\left\{B\left(x,\tfrac{1}{2^{N_\varepsilon}}\right)\:\middle|\: x\in\llbracket T\rrbracket\right\}=\left\{\llbracket T_p\rrbracket\:|\: p\in T\text{ and }\len(p)=K_\varepsilon\right\}
	\end{equation*}
	must be finite. Therefore $(\llbracket T\rrbracket,d_2)$ is totally bounded, and since it is also complete (by Proposition \ref{TreeTopologyUltrametrizable}) it is compact.\\
	Suppose now that $\left(\llbracket T\rrbracket,d_2\right)$ is compact, and thus totally bounded. Suppose by contradiction that there exists a position $p_{_0}\in T$ such that $p_{_0}$ has an infinite number of offsprings. Denote $K_0=\len(p_{_0})$. There cannot exists a finite number of open balls of radius smaller than $\tfrac{1}{2^{K_0}}$ such that their union covers $\left\llbracket T_{p_{_0}}\right\rrbracket$ since it is equal to a union of infinitely many open balls of radius $\tfrac{1}{2^{K_0}}$, a contradiction to the assumption that $\left(\llbracket T\rrbracket,d_2\right)$ is totally bounded.
\end{proof}
\begin{corollary}
	\label{TotallyBoundedIffFiniteCorollary}
	Let $(X,\mathcal{O})$ be an ultrametrizable space, and let $T_X$ be a representing tree of $(X,\mathcal{O})$. Then $(X,\mathcal{O})$ is totally bounded if and only if every position of $T_X$ has a finite number of offsprings.
\end{corollary}
\begin{proof}
	Since an ultrametrizable space and its completion share the same representing trees, by Proposition \ref{CompactTree} $(X,\mathcal{O})$ is totally bounded if and only if its completion is compact if and only if its representing trees have positions with finite numbers of offsprings.
\end{proof}
\begin{remark}
	The fact that a totally bounded space is also separable is an immediate consequence of Corollary \ref{TotallyBoundedIffFiniteCorollary}, in the case of ultrametric spaces.
\end{remark}
\begin{corollary}
	Let $T$ be a pruned tree, and let $\emptyset\neq P\subseteq\llbracket T\rrbracket$. The clade space $(P,\widetilde{\mathcal{O}}_T)$ is locally totally bounded if and only if every branch $x\in P$ has a finite prefix $p_x\in T$ such that every position in the subtree $T_{p_x}$ has a finite number of offsprings.
\end{corollary}

\subsection{Discrete topology and perfectness}
A topological space is discrete if every point in it is an isolated point. A topological space is perfect if it has no isolated points. The identification of isolated points in a clade space are simple:
\begin{proposition}
	Let $T$ be a tree and let $x\in\llbracket T\rrbracket$. Then $x$ is an isolated point (according to the tree topology) if and only if all by finitely many of its prefixes have only one offspring. 
\end{proposition}
\begin{proof}
	Let $T$ be a tree and let $x\in\llbracket T\rrbracket$. Let $d_2$ be the dyadic ultrametric over the canopy $\llbracket T\rrbracket$. Recall that $d_2$ generates the tree topology.\\
	Then $x$ is an isolated point if and only if there exists a $K\in\Naturals$ such that the open balls $B\left(x,\tfrac{1}{2^{k}}\right)=\{x\}$ for every $k\geq K$, which is true if and only if all prefixes off length at least $K$ have only one offspring.
\end{proof}
\begin{corollary}
	A clade space is discrete if and only if almost every finite prefixes of every of its branches have only one offspring. A clade space is perfect if it has no branch such that almost every of its finite prefixes have only one offspring.
\end{corollary}

\subsection{Doubling sequentially-descending-to-zero ultrametric}
The property of being \emph{doubling} has a simple characterization in the case of a sequentially-descending-to-zero ultrametric, when using the representing tree.
\begin{proposition}
	\label{DoublingImpliesBoundedOffsprings}
	Let $(X,d)$ be a sequentially-descending-to-zero ultrametric space, and let $T_{(X,d)}$ and $\left(\llceil X\rrceil,d_{T_{(X,d)}}\right)$ be its representing tree and clade space respectively. If $(X,d)$ is doubling with constant $D\in\Naturals\setminus\{0\}$, then every position of $T_{(X,d)}$ has no more than $D$ offsprings.
\end{proposition}
\begin{proof}
	Suppose that $(X,d)$ is doubling with constant $D\in\Naturals\setminus\{0\}$. Thus by Proposition \ref{IsometricToTreeTheorem} also $\left(\llceil X\rrceil,d_{T_{(X,d)}}\right)$ is doubling with the constant $D$. Suppose by contradiction that its representing tree $T_{(X,d)}$ has a position $p^*$ with more than $D$ offsprings. Denote $k^*=\len(p^*)$. Let $(r_k)_{k\in\Naturals}$ such as described in Definition \ref{representingDefinition} and let $x\in\llceil X\rrceil$ such that $p^*\prec x$. Then
	\begin{equation}
		\label{CoveringDendronicBall}
		\closedb(x,r_{k^*})=\bigcup_{y\in\closedb(x,r_{k^*})}\closedb(y,r_{k^*+1}),
	\end{equation}
	and since $p^*$ has more than $D$ offsprings, the union in the right hand side of Equation (\ref{CoveringDendronicBall}) consists of more than $D$ non-intersecting distinct closed balls. Denote $r=\tfrac{r_{k^*}+r_{k^*+1}}{2}$. Then
	\begin{equation*}
		\closedb(y,r)=\closedb(y,r_{k^*+1})
	\end{equation*}
	for every $y\in\llceil X\rrceil$, and
	\begin{equation*}
		\closedb(x,2r)\supseteq\closedb(x,r_{k^*}).
	\end{equation*}
	Thus, the closed ball $\closedb(x,2r)$ cannot be covered by no more than $D$ closed ball of radius $r$, a contradiction. Hence every position of $T_{(X,d)}$ has no more than $D$ offsprings.
\end{proof}
\begin{corollary}
	Let $(X,d)$ be a doubling ultrametric space. Then $(X,d)$ is totally bounded.
\end{corollary}
\begin{proof}
	This is a consequence of Corollary \ref{TotallyBoundedIffFiniteCorollary} and Proposition \ref{DoublingImpliesBoundedOffsprings}.
\end{proof}
While Proposition \ref{DoublingImpliesBoundedOffsprings} gives a necessary condition for a sequentially-descending-to-zero ultrametric to be doubling, Proposition \ref{DoublingDendronicProposition} gives some sufficient conditions for a sequentially-descending-to-zero ultrametric spaces to be doubling, based on the image of the sequentially-descending-to-zero ultrametric and on its representing tree.
\begin{proposition}
	\label{DoublingDendronicProposition}
	Let $(X,d)$ be a sequentially-descending-to-zero ultrametric space and let $(r_k)_{k\in\Naturals}$ be a weakly monotonic decreasing sequence of non-negative numbers such that
	\begin{equation*}
		\{0\}\cup\{r_k \:|\: k\in\Naturals\}=\image d
	\end{equation*}
    and $r_k>r_{k+1}$ if $r_k>0$, for every $k\in\Naturals$.\\
    Additionally, let $T_{(X,d)}$ and $\left(\llceil X\rrceil,d_{T_{(X,d)}}\right)$ be the representing tree and clade space of $(X,d)$, respectively. Then $(X,d)$ is doubling if there is a finite upper bound on the number of offsprings of each position in $T_{(X,d)}$ and at least one of the following holds:
	\begin{enumerate}
		\item\label{eventuallyOneOffspring} There is a length $K_1\in\Naturals$ such that all positions $p\in T_{(X,d)}$ of length $\len(p)\geq k_1$ have only one offspring.
		\item\label{ratioBounded} The sequence $(r_k)_{k\in\Naturals}$ has only positive elements and there exists an $l\in\Naturals$ such that for every $k\in\Naturals$
		\begin{equation*}
			\tfrac{r_k}{r_{k+k}}\geq 2.
		\end{equation*}
	\end{enumerate}
\end{proposition}
\begin{proof}
	Suppose that $(X,d)$ satisfies the condition from point \ref{eventuallyOneOffspring}. Since the number of offspring is finite for every position, define
	\begin{equation*}
		\begin{split}
			&M= \max\left\{|\{p'\:|\: p\prec p'\text{ and }\len(p')=\len(p)+1\}| \: \middle|\: p\in T_{(X,d)}\right\} \\
			&=\max\left\{|\{p'\:|\: p\prec p'\text{ and }\len(p')=\len(p)+1\}| \: \middle|\: p\in T_{(X,d)}\text{ and }\len(p)<k_1\right\}.\\
		\end{split}
	\end{equation*}
	Then, $\left(\llceil X\rrceil,d_{T_{(X,d)}}\right)$ is doubling, since it has at most $M\cdot k_1$ elements, and every finite metric space is doubling.\\
	Suppose now that $(X,d)$ satisfies the conditions from point \ref{ratioBounded}.
	Denote by $M\in\Naturals\setminus\{0\}$ the upper bound on the number of offsprings of each position of $T_{(X,d)}$.	We claim that $\left(\llceil X\rrceil,d_{T_{(X,d)}}\right)$ is doubling. Indeed, let $x\in\llceil X\rrceil$ and $r\in(0,\infty)$. Denote
	\begin{equation*}
		K_1=\max\{k\in\Naturals\:|\: r_k\leq r\}
	\end{equation*}
	and
	\begin{equation*}
		K_2=\max\{k\in\Naturals\:|\: r_k\leq 2r\}.
	\end{equation*}
	If $K_1\leq l$, then there are no more than $M^l$ positions of length $K_1$. Thus also the position 
	\begin{equation*}
		\closedb(x,r_{K_2})=\closedb(x,2r)
	\end{equation*}
	can be covered with no more than $M^l$ closed balls of radius $r$ (or equivalently $r_{K_1}$).
	Else, suppose that $K_1\geq 1$. We have 
	\begin{equation*}
		2r_{K_1-1-l}\geq 2r_{K_1-1}\geq2r,
	\end{equation*}
	and thus
	\begin{equation*}
		r_{K_1-1-l}\geq r_{K_2}.
	\end{equation*}
	Hence,
	\begin{equation*}
		K_1-1-l\leq K_2,
	\end{equation*}
	so that
	\begin{equation*}
		K_1-K_2\leq 1+l.
	\end{equation*}
	Thus the position $\closedb(x,2r)=\closedb(x,r_{K_2})$ has no more than $M^{K_1-K_2}$ descendants of length $K_1$, which is no more than $M^{l+1}$ descendants of that length. Thus the closed ball $\closedb(x,2r)$ can be covered by no more than $M^{l+1}$ closed balls of radius $r$. In both cases, $\left(\llceil X\rrceil,d_{T_{(X,d)}}\right)$ is indeed doubling with the doubling constant $M^{l+1}$, and by Theorem \ref{IsometricToTreeTheorem}, so is $(X,d)$.
\end{proof}
    
%\printbibliography
\bibliographystyle{plain}
\bibliography{masterthesis.bib}

\begin{thebibliography}{1}

\bibitem{bellaiche2023}
I.~Bella{\"i}che.
\newblock Dyadic hausdorff dimension games.
\newblock Master's thesis, Tel Aviv University, 2023.

\bibitem{brian2015completely}
W.R. Brian.
\newblock Completely ultrametrizable spaces and continuous bijections.
\newblock In {\em Topology Proceedings}, volume~45, pages 233--252, 2015.

\bibitem{dovgoshey2022range}
O.~Dovgoshey and V.~Scherbak.
\newblock The range of ultrametrics, compactness, and separability.
\newblock {\em Topology and its Applications}, 305:107899, 2022.

\bibitem{hughes2004trees}
B.~Hughes.
\newblock Trees and ultrametric spaces: a categorical equivalence.
\newblock {\em Advances in mathematics}, 189(1):148--191, 2004.

\bibitem{mauldin2009graph}
R.D. Mauldin, T.~Szarek, and M.~Urba{\'n}ski.
\newblock Graph directed markov systems on hilbert spaces.
\newblock In {\em Mathematical Proceedings of the Cambridge Philosophical Society}, volume 147.2, pages 455--488. Cambridge University Press, 2009.

\end{thebibliography}
\end{document}